\documentclass[11pt]{article}
\usepackage[backrefs, alphabetic,initials]{amsrefs}
\usepackage{amsfonts,amsmath,amscd,amssymb,amsthm}
\usepackage{mathrsfs}
\def\0{{\bar 0}}
\def\1{{\bar 1}}

\def\ch{{\operatorname{ch}\:}}

\def\lm{{\operatorname{LM}}}
\def\LM{{\operatorname{\bf LM}}}

\newcommand{\noi}{\noindent}
\newcommand{\ga}{\alpha}

\newcommand{\gc}{\gamma}
\newcommand{\gd}{\delta}

\newcommand{\gs}{\sigma}

\newcommand{\gl}{\lambda}

\newcommand{\gr}{\rho}

\newcommand{\gep}{\epsilon}

\def\RR{{\mathtt{R}}}
\def\SS{{\mathtt{S}}}

\newcommand{\ci}{\circ} \newcommand{\ti}{\times}
\newcommand{\fg}{\FRAK{g}}\newcommand{\fgl}{\FRAK{gl}}
\newcommand{\fsl}{\FRAK{sl}}

\newcommand{\fh}{\FRAK{h}}

\newcommand{\fb}{\FRAK{b}}

\newcommand{\fp}{\FRAK{p}}
\newcommand{\fq}{\FRAK{q}}
\newcommand{\fl}{\FRAK{l}}

\newfont{\eufm}{eufm10 scaled\magstep1}
\newcommand{\FRAK}[1]{\mbox{\eufm#1}}

\newcommand{\cP}{\mathcal{P}}
\newcommand{\cF}{\mathcal{F}}

\newcommand{\cG}{\mathcal{G}}

\newcommand{\cR}{\mathcal{R}}

\newcommand{\bco}{\begin{conjecture}}
\newcommand{\ba}{\begin{alg}}
\newcommand{\ea}{\end{alg}}
\newcommand{\eco}{\end{conjecture}}
\newcommand{\bpf}{\begin{proof}}
\newcommand{\epf}{\end{proof}}
\newcommand{\bt}{\begin{theorem}}
\newcommand{\et}{\end{theorem}}

\newcommand{\br}{\begin{rem}}
\newcommand{\er}{\end{rem}}
\newcommand{\bi}{\begin{itemize}}
\newcommand{\bl}{\begin{lemma}}
\newcommand{\bsul}{\begin{sublemma}}
\newcommand{\esul}{\end{sublemma}}
\newcommand{\bp}{\begin{proposition}}
\newcommand{\be}{\begin{equation}}
\newcommand{\bc}{\begin{corollary}}
\newcommand{\bexa}{\begin{example}}
\newcommand{\eexa}{\end{example}}
\newcommand{\bex}{\begin{exercise}}
\newcommand{\eex}{\end{exercise}}
\newcommand{\btab}{\begin{tab}}
\newcommand{\etab}{\end{tab}}
\newcommand{\brs}{\begin{rems}}
\newcommand{\ers}{\end{rems}}
\newcommand{\ei}{\end{itemize}}
\newcommand{\el}{\end{lemma}}
\newcommand{\ep}{\end{proposition}}
\newcommand{\ee}{\end{equation}}
\newcommand{\ec}{\end{corollary}}
\newcommand{\Bc}{\begin{center}}
\newcommand{\Ec}{\end{center}}

\newcommand{\sfl}{\small{\mathfrak{l}}}
\newcommand{\sfp}{\small{\mathfrak{p}}}
\newcommand{\sfq}{\small{\mathfrak{q}}}
\newcommand{\sfb}{\small{\mathfrak{b}}}

\newcommand{\sfg}{\small{\mathfrak{g}}}

\def\Ext{{\operatorname{Ext}}}
\numberwithin{equation}{section}%
\newcommand{\lra}{\longrightarrow}

\begin{document}
\title{Table of Contents}
\newcommand{\rh}{\epsilon^{\widehat{\rho}}}
\newcommand{\pii}{\prod_{n=1}^\infty}
\newcommand{\sii}{\sum_{n= - \infty}^\infty}
\newtheorem{theorem}{Theorem}[section]
 \newtheorem{lemma}[theorem]{Lemma}
 \newtheorem{example}[theorem]{Example}
  \newtheorem{sublemma}[theorem]{Sublemma}
   \newtheorem{corollary}[theorem]{Corollary}
 \newtheorem{conjecture}[theorem]{Conjecture}
 \newtheorem{hyp}[theorem]{Hypothesis}
 \newtheorem{alg}[theorem]{Algorithm}
  \newtheorem{sico}[theorem]{Sign Convention}
 \newtheorem{tab}[theorem]{\quad \quad \quad \quad \quad \quad \quad \quad \quad \quad \quad \quad \quad Table }
  \newtheorem{rem}[theorem]{Remark}
 \newtheorem{proposition}[theorem]{Proposition}
\font\twelveeufm=eufm10 scaled\magstep1 \font\teneufm=eufm10
\font\nineeufm=eufm9 \font\eighteufm=eufm8 \font\seveneufm=eufm7
\font\sixeufm=eufm6 \font\fiveeufm=eufm5
\newfam\eufmfam
\textfont\eufmfam=\twelveeufm \scriptfont\eufmfam=\nineeufm
\scriptscriptfont\eufmfam=\sixeufm \textfont\eufmfam=\teneufm
\scriptfont\eufmfam=\seveneufm \scriptscriptfont\eufmfam\fiveeufm
\newtheorem{exercise}{}
\setcounter{exercise}{0} \numberwithin{exercise}{section}
\newenvironment{emphit}{\begin{theorem}}{\end{theorem}}
\setcounter{Thm}{0} \numberwithin{Thm}{section} 
\title{\large Combinatorics of Character
Formulas for the Lie Superalgebra $\fgl(m,n).$}
\author{Ian M. Musson
\\Department of Mathematical Sciences\\
University of Wisconsin--Milwaukee\\ and\\
Vera V. Serganova\\
Department of Mathematics\\
University of California--Berkeley.\\} \maketitle
\begin{abstract} Let $\fg$ be the Lie superalgebra $\fgl(m,n).$  Algorithms for computing the
composition factors and multiplicities  of Kac modules for $\fg$
 were given  by the second author, \cite{S2}  and by J. Brundan \cite{Br}.

We give a combinatorial proof of the equivalence between the two
algorithms. The proof uses weight and cap diagrams introduced by
Brundan and C. Stroppel, and cancelations between paths in a graph $\mathcal{G}$ defined using these diagrams.
Each vertex of $\mathcal{G}$ corresponds to a highest weight of a finite dimensional simple module, and each edge  is weighted by a nonnegative integer. If $\mathcal{E}$ is the subgraph of  $\mathcal{G}$ obtained by   deleting all edges of positive weight, then $\mathcal{E}$ is the graph that describes non-split extensions between simple highest weight modules.

We also give a procedure for finding the composition factors of
any Kac module, without cancelation. This procedure leads to a second proof of the main result.\end{abstract}

\section{Introduction.}
The problem of finding the characters of the finite  dimensional
simple modules for the complex Lie superalgebra $\fg = \fgl(m,n)$
was first posed by V. Kac in 1977, \cite{Kac1}.  Let $X^+(m,n)$ denote the set of dominant
integral weights for $\fg$.  In  \cite{Kac2} Kac introduced a certain finite dimensional
highest weight module $K(\gl)$, now known as a
Kac module, with highest weight   $\gl \in X^+(m,n)$, whose character is given by an analog of  the Weyl character formula.  Furthermore any composition factor of $K(\gl)$ is a simple module $L(\mu)$ with highest weight $\mu \in X^+(m,n),$ and the multiplicities of the composition factors of Kac modules can be expressed using an upper triangular matrix with diagonal entries equal to 1.   Therefore the determination of this multiplicity
matrix leads to a solution of the problem raised by Kac.

Combinatorial formulas for the
multiplicity of $L(\lambda)$  as a composition factor of $K(\mu)$
 were given  in  \cite{S2}  and
\cite{Br}, using completely different methods.  We give a combinatorial proof of the equivalence between these two formulas, Theorem A.
Let $F$ be
the set of all functions from $\mathbb{Z}$
 to the set $ \{ \ti, \ci, < \;,\; > \}$ such that $f(a) = \ci$ for all but finitely
 many $a \in \mathbb{Z}$.
Let $\mathbb{Z}F$ be the free abelian group with basis $F$.
The idea of the proof is to express the formula from
\cite{S2}  as a signed sum of terms in $\mathbb{Z}F$.
The terms from this sum correspond to
paths in a certain graph $\mathcal{G}$ .  We define an involution on the paths occurring in this sum such that paths that are paired by the involution have opposite signs.   After these terms are canceled, what remains is the formula from \cite{Br} in a form communicated to the first author by Brundan.  This reformulation of  the result from \cite{Br} uses diagrams called weight and cap diagrams that originate in the work of Brundan and Stroppel on Khovanov's diagram algebra, \cite{BS1}, \cite{BS2}, \cite{BS3}, \cite{BS4}.   We remark that our notation for these diagrams is different from theirs.
We note also that character formulas for the irreducible representations of the orthosymplectic Lie superalgebras were announced in \cite{S3}.  These results are expressed in terms of weight diagrams and proved in \cite{GS}.
\\ \\
Since the category of finite dimensional $\mathbb{Z}_2$-graded weight modules $\cF$ is not semisimple, an important problem in representation theory is to determine the non trivial extensions between simple modules.  This problem is related to the graph $\cG$ as follows.
Each vertex $f$ of $\mathcal{G}$ corresponds to a highest weight of a finite dimensional simple module $L(f)$, and each edge of $\mathcal{G}$ is weighted by a nonnegative integer. In Theorem B we show that if $\mathcal{E}$ is the subgraph of  $\mathcal{G}$ obtained by deleting all edges of positive weight, then $\Ext^1_\cF(L(f),L(g)) \neq 0$ if and only
if $f\longrightarrow g$ or $g\longrightarrow f$ is an edge of $\mathcal{E}$.

This paper is organized as follows.  In the next section we give a formal statement of the main combinatorial result, Theorem A.  Some work is necessary to derive the equivalence of the character formulas from the combinatorial statement, and this is done in Section 3.  The graph $\cG$ is introduced and Theorem A is proved in Section 4.
In Section 5 we outline a procedure for finding the composition factors of
any Kac module, without cancelation. This procedure leads to a second proof of the main result.  Theorem B relating extensions to the subgraph $\mathcal{E}$ is presented in Section 6.

Although some of our results can be deduced from work of Brundan and Stroppel,  we have included our original proofs, so that our work may
be read independently.  For further details, see Remark \ref{rm69}.

Our notation for weight diagrams is the same as in \cite{GS} and
hence slightly different from those in \cite{BS4}.
Here is a ``dictionary'' which allows quickly pass from one set of notation
to the other. \\
\\
Our notation $\quad \quad \quad \quad \quad \quad \quad \quad \times \quad \quad \quad \circ \quad \quad \quad > \quad \quad \quad < $
\\ \\
Brundan-Stroppel  notation  $\quad \quad \vee \quad \quad \quad \wedge \quad \quad \quad  \times \quad \quad \quad \circ$
\\
\\
The authors would like to express their gratitude to Jon Brundan for sharing his ideas with them, and to thank Hiroyuki Yamane for pointing out an error in an earlier version of this paper.

\section{A Combinatorial Formula.}
 \indent Let $F$ be
the set of all functions from $\mathbb{Z}$
 to the set $ \{ \ti, \ci, < \;,\; > \}$ such that $f(a) = \ci$ for all but finitely
 many $a \in \mathbb{Z}$. 
 For $f \in F$ we set
 $\#f = |f^{-1}(\ti)|,$  and \[ core_L(f) = f^{-1}(>), \quad core_R(f) = f^{-1}(<).\]
We call $\#f$ the {\it degree of atypicality} of $f$, and define the
the {\it core} of $f$ to be  $$core(f) = (core_L(f),core_R(f)).$$
Let $\mathbb{Z}F$ be the free abelian group with basis $F$.  Our
main result is an identity for certain $\mathbb{Z}$-linear operators
on $\mathbb{Z}F$. If  $f \in F$ we define the  {\it weight diagram}
$D_{wt}(f)$ to be a  number line with the symbol $f(a)$ drawn at
each $a \in \mathbb{Z}$. Next let $C_L$ and $ C_R$ be disjoint
finite subsets of $\mathbb{Z}$ and consider a number line with
symbols $>$ (resp. $<$) located at all $a \in C_L$ (resp. $a \in
C_R)$. A {\it cap} $C$ is the upper half of a circle joining two
integers $a$ and $b$ which are not in $C_L \cup C_R$.  If $b < a$ we
say that $C$ {\it begins} at $b$ and {\it ends} at $a$ and we write
$b(C) = b$, and $e(C) = a.$ A finite set of caps, together with the
symbols $< \;, \;>$ located as above, is called a {\it cap diagram}
if no two caps intersect, and the only integers inside the caps, which
are not ends of some other caps, are
located at points in $C_L \cup C_R$.  We remark that in \cite{BS4} the last condition is ensured as follows: infinite vertical rays are drawn at all vertices that are not in $C_L \cup C_R$ and not ends of a cap.  Then the requirement is that no two distinct caps or rays can intersect.

\indent If $D$ is a cap diagram there is a unique $f \in F$ such that $core(f) = (C_L,C_R)$, and for $a \notin C_L \cup C_R$,
\begin{eqnarray*}
f(a) & = & \ti
\; \mbox{if there is a cap in $D$ beginning at} \; a,\\
f(a) & = & \ci \; \; \mbox{otherwise}.
\end{eqnarray*}

\indent We write $D = D_{cap}(f)$ in this situation.
If $D = D_{cap}(f)$ or $D = D_{wt}(f)$ we set $core(D) = core(f).$\\
\\
We say that a weight diagram and a cap diagram {\it match} if they
have the same core and, when superimposed on the same number line,
each cap connects a $\ti$ to a $\ci$.  For $f \in F$, set \be
\label{Pfdef}   P(f) = \{g \in F|D_{cap}(g) \; \mbox{matches} \;
D_{wt}(f)\} . \ee Brundan's formula for the composition factors of a
Kac module can be written in terms of matching cap and weight
diagrams. We now turn to the combinatorics necessary to express
 the formula from \cite{S2}. If $f \in F$ and $\# f = k$ we set \be
\label{e2} \ti(f) = (a_1, a_2, \ldots,a_k) \ee if $f^{-1}(\ti) =
\{a_1, a_2, \ldots, a_k\}$ with $a_1 >a_2 > \ldots >a_k$.\\
\\
Next suppose  $f$ satisfies (\ref{e2}), and that  $f(a) = \ti,$ and
$f(b)  =  \ci.$ Informally, we define $f_{b} \in F,$ (resp. $f^a \in
F$) by adding $b$ to $\ti(f)$ (resp. deleting $a$ from $\ti(f)$).
Precisely $f_{b}$ and $f^a$ have the same core as $f$, and satisfy
\[\ti(f_{b}) =  (a_1, \ldots, a_j, b, a_{j+1}, \ldots, a_k), \]
\[\ti(f^a) =  (a_1, \ldots, a_{i-1},a_{i+1},\ldots, a_k),\]
where $a = a_i,$ and $a_j > b > a_{j+1}.$ Here it is convenient to
set $a_0 = \infty,$ and $a_{k+1} =- \infty.$ We also set $f_{b}^a =
(f_{b})^a =  (f^a)_{b}.$\\
\\
For $f \in F,$ and $a, b \in \mathbb{Z}$ with $b<a,$ let
$l_f(b,a)$ be the number of occurrences of the symbol
$\ti$ minus the number of occurrences of $\ci$ strictly between $b$
and $a$ in the  weight diagram of $f$. We say that  $g \in F$ is
obtained from $f$ by a {\it legal move} if $g = f^a_b$ for some $b<a$
such that $f(a) =\ti, f(b)=\circ$ and $l_f(b,c) \geq 0$ for all $c$ with $b< c \leq a$. We call $a$ the {\it start}, $b$ the {\it end} and $l_f(b,a)$ the {\it weight} of the legal move.\\
\\
There is another way to think about legal moves. Suppose that
 $b < a,$  $f(b) = \ci$ and $f(a) = \ti$.  Keep a tally starting at $b$ with a tally of zero, and move to the right  along the number line adding one to the tally every time a $\ti$ is passed, and subtracting one every time a $\ci$ is passed in the  weight diagram
$D_{wt}(f)$ . Then $f_{b}^a$ is
obtained from $f$ by a legal move  if  and only if the tally remains non-negative on the interval $[b,a)$.  If  this is the case the weight of the  legal move is the value of the tally just before we arrive at $a.$

\begin{example} {\rm Assume that $f$ is core-free with $\times(f) = (1, 2,4,5,7)$ and let $g = f_0^7.$  There is a legal move from $f$ to $g$ of weight 2. We reflect the cap diagram $D_{cap}(g)$ in the number line,
thus obtaining a cup diagram, which is drawn in place of $D_{cap}(g).$  Then we superimpose the usual cap diagram 
$D_{cap}(f)$.  The resulting diagram is similar to the diagrams of Brundan and Stroppel. A significant difference is that our diagram is not  oriented.}
$$  
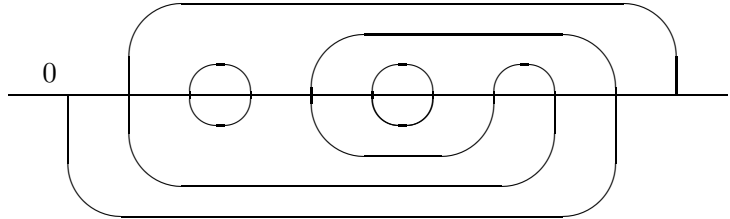
\begin{picture}(8,105)
\put(-120.8,50){$0$}
\put(-133.8,45){\line(1,0){272}}
\put(-7.5,45){\oval(161,69)[b]}
\put(38.5,45){\oval(115,46)[t]}
\put(61.5,45){\oval(23,23)[t]}
\put(-53.5,45){\oval(23,23)}
\put(15.5,45){\oval(23,23)}
\put(15.5,45){\oval(23,23)[b]}
\put(15.5,45){\oval(207,69)[t]}
\put(-7.5,45){\oval(207,92)[b]}
\put(15.5,45){\oval(69,46)[b]}
\end{picture}$$
\vspace{0.1cm}
\end{example}
It follows immediately from the definition of a legal move that $g$ is
obtained from $f$ by a legal move of weight $0$ (starting at $a$ and
ending at $b$) if and only if $D_{cap}(g)$ has a cap joining vertex
$a$ and vertex $b$, where $g(b) =\ti$ and $g(a) =\ci$, and $g$ is obtained from $f$ by switching the labels at
these two vertices.
\\ \\
Next if $1 \leq i \leq k$ we define operators
$\sigma_i:\mathbb{Z}F \longrightarrow \mathbb{Z}F $ as follows. If
$f$ satisfies (\ref{e2}) and $a = a_i,$ then \be \label{sif}
\sigma_i(f) = \sum_b (-1)^{l_f(b,a)} f_{b}^a
\;,\ee where the sum is over all $b$ such that $f_{b}^a$ is obtained from $f$ by a legal move.\\
\\
 Now we can state our first main result.\\ \\
{\bf Theorem A.} \label{T1} For $f \in F$ with \#$f=k$ we have \be \label{eq 1}
(1 + \sigma_1) \ldots (1 + \sigma_k)f = \sum_{g \in P(f)} g. \ee

\bexa {\rm Suppose that $core(f) =
(\{7\}, \{3\})$ and
\[ \ti(f) = (9, 6, 5, 1, 0), \]
then $D_{cap}(f)$ is pictured below.  Note that $k = 5$ and $a_1 =
9.$

$$
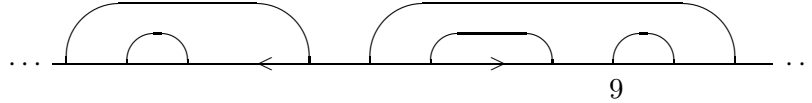
\begin{picture}(130,24)
\put(-91,0.5){$\cdots$} \put (-74,3){\line(1,0){272}}
\put(203,0.5){$\cdots$} \put(-23,3){\oval(92,46)[t]}
\put(-34.5,3){\oval(23,23)[t]} \put(115,3){\oval(138,46)[t]}
\put(149.5,3){\oval(23,23)[t]} \put(92,3){\oval(46,23)[t]}
\put(3,1){$\scriptstyle<$} \put(91,1){$\scriptstyle>$}
\put(136.5,-10){9}
\end{picture}
$$

\vspace{0.5cm}

Now $f_{b}^9$ is
obtained from $f$ by a legal move if and only if $b = -1,4$ or 8.  The weights of these legal moves are 1, 1 and 0 respectively. Thus Equation
(\ref{sif}) becomes
$$\gs_1(f) = f^9_{8} - f^9_{4} - f^9_{-1}.$$

Replacing 9 by $-1, 4, 8$ in $\ti(f)$ we obtain 
\[ \ti(f^9_{-1}) =(6, 5, 1, 0, -1),
\] \[\ti(f^9_4) = (6, 5, 4, 1, 0),\]
\[ \ti(f^9_8) = (8, 6, 5, 1, 0) .\]

The cap
diagram $D_{cap}(f_4^9)$ is given below. 
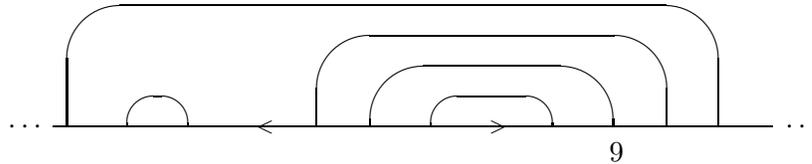
\vspace{1.0cm}
$$
\begin{picture}(130,24)
\put(-91,0.5){$\cdots$} \put(-74,3){\line(1,0){272}}
\put(203,0.5){$\cdots$} \put(54.5,3){\oval(246.5,92)[t]}
\put(-34.5,3){\oval(23,23)[t]}\put(92,3){\oval(132.5,69)[t]}
\put(92,3){\oval(46,23)[t]} \put(3,1){$\scriptstyle<$}
\put(91,1){$\scriptstyle>$} \put(136.5,-10){9}
\put(92,3){\oval(92,46)[t]}
\end{picture}
$$ }
\vspace{0.1cm} 

\eexa \br \label{rem2.1} {\rm The real content of Theorem A is the core-free
case.  Indeed since $\sigma_i f$ is a linear combination of terms
$h \in F$ with the same core as $f$, and any $g \in P(f)$, has the
same core as $f$, we immediately reduce to this case.  The symbols $<$ and $>$  are important in the application to the Lie superalgebra $\fgl(m,n).$
When $f$ is core-free then $f$
is completely determined by $\ti(f)$}. \er 
\bexa \label{catalan}{\rm An interesting case arises when $f$ is
  core-free and $$\ti(f) =
(2,4,6,\ldots,2k-2).$$ Let $A_k$ be the set of cap diagrams with $k$
caps each of which begins and ends at points in the set $ \{0,1,2,
\ldots, 2k-1\},$ and let $B_k$ be the set of cap diagrams that match
the weight diagram $D_{wt}(f).$ Given a diagram in $A_k,$ we obtain
a diagram in $B_k$ by deleting the cap beginning at 0. This gives a
bijection from $A_k$ to $B_k$. The cardinality of $A_k$ is the $kth$
Catalan number $C_k = \frac{1}{k+1} \left(\begin{array}{c} 2k\\ k
\end{array} \right) $ see \cite{stan} Exercise 6.19 part o. In terms of representation theory this
means that if $\gr$ is defined as in (\ref{nrd}) below, then  the  length
of a composition series for the Kac module $K(\gr)$ for $\fgl(k-1,k-1)$ equals $C_k$.  There are further examples in  \cite{su} where the number of
composition factors of a Kac module is a Catalan number.  We conjecture that if $g$ is core-free and $|g^{-1}(\ti)| = k-1$, then $|P(g)| \leq C_k$ with equality if and only if $\ti(g)$ is obtained from $\ti(f)$ by adding the same integer to each entry.}\eexa
\section{Character Formulas.}

\indent In this subsection $\fg$  will be the Lie superalgebra
$\mathfrak{gl}(m,n)$ and
 $\fh$
and $\fb$ the  Cartan and Borel subalgebras, consisting of diagonal
and upper triangular matrices respectively. Let $\epsilon_{i}, \,
\delta_{j}$ be the linear functionals on $\FRAK{h}$ whose value on
the diagonal matrix
\[ a = diag(a_{1}, \ldots, a_{m+n}) \]
is given by \be \label{eddef}  \epsilon_{i}(a) = a_{i}, \;\;
\delta_{j}(a) = a_{m+j} \;\;
 1 \leq i \leq m, \; 1 \leq j \leq n. \ee
We define a bilinear form $(\;,\;)$ on $\fh^*$ by \be
\label{edform}(\epsilon_i,\epsilon_j) = \delta_{i,j} = - (\delta_i,
\delta_j).\ee
Let $X = X(m|n)$ denote the lattice of integral
weights spanned by the $\epsilon_i$ and $\delta_i.$
Also set
\be \label{typea}  \Delta_{0}^+ = \{ \epsilon_{i} - \epsilon_{j};
\delta_{i} - \delta_{j} \}_{i < j}, \; \Delta_{1}^+ = \{
\epsilon_{i} - \delta_{j} \} \; \mbox{and}\; \Delta^+ = \Delta^+_{0} \cup\Delta^+_{1}. \ee
Then $\Delta^+$ is the set of roots of $\fb.$  Next let
\be \label{nrd} \rho = m \epsilon_1
+\cdots+2\epsilon_{m-1}+\epsilon_m-\delta_1-2\delta_2-\cdots-n\delta_n.
\ee A weight $\gl
\in X$ is {\it regular} if $(\gl+\rho,\epsilon_i - \epsilon_j) \neq 0,$ and
$(\gl+\rho, \delta_i - \delta_j)\neq 0$ if $i \neq j.$ Let $X_{reg}$ be
the subset of $X$ consisting of regular weights, and
let $S_m$ be the symmetric group of degree $m$. The Weyl group
$W = S_m \times S_n$ acts on $\fh^*$ by permuting the $\gep_i$ and
$\gd_i.$ The dot action of $W$ is defined by $w\cdot \gl =
w(\gl+\gr) - \gr.$ We will identify $\gl \in X(m|n)$ with the $m+n$
tuple of integers
\be \label{alb} (a_1 , a_2 , \cdots , a_m| b_1 , b_2 , \cdots , b_n)\ee where \be
\label{lab} a_1 = (\lambda+\rho,\epsilon_1), \cdots, a_m =
(\lambda+\rho,\epsilon_m),\quad b_1 = (\lambda+\rho,\delta_1) ,
\cdots , b_n =  (\lambda+\rho,\delta_n). \ee
Let $X^+ = X^+(m|n)$
denote the set of $\gl = (a_1 , a_2 , \cdots , a_m| b_1 , b_2 ,
\cdots , b_n)$ in $X(m|n)$ with \be \label{X0}  a_1 > a_2 > \cdots > a_m, \quad
b_1 < b_2 < \cdots < b_n.\ee  In this $\gr$-shifted notation, the dot
action of $W$ is represented by permutations of the entries in
$\gl.$ If $\gl \in X_{reg}$ there is a unique element $w$ of the Weyl group
$W$ such that $w\cdot \gl \in X^+.$\\
\\
   Given $f \in F,$ denote by $\gl = \gl(f)$ the element of $X^+$ such that when written in the form (\ref{alb}), the entries of on the left (resp. right) side of $\gl$ coincide with $ core_L(f) \cup \ti(f)$  (resp. $core_R(f) \cup \ti(f)$)   arranged  in order as in (\ref{X0}).
    This defines a bijection
$\gl \lra f_\gl$
   from $X^+$ to $F$ whose inverse we write as $f \lra \gl(f).$
\noi Let $\cF$ be the category of finite
dimensional $\fg$ modules which are weight modules for $\fh,$ and
for $\gl \in X^+,$ let $K(\gl)$ (resp. $L(\gl)$) be the Kac module
(resp. simple module) with highest weight $\gl$.  The map $f \lra L(\gl(f))$     extends to an isomorphism from  $\mathbb Z F$ to the Grothendieck group of $\cF$, and we often identify these two groups.
 For $\gl \in X^+,$  we write $D_{wt}(\gl), \;D_{cap}(\gl),$ $ \# \gl$ and   $\times(\gl)$
in place of $D_{wt}(f_\gl), \;D_{cap}(f_\gl),$ $ \#(f_\gl)$ and
$\times(f_\gl)$ respectively, and set
\be \label{pmgl}  \mathbf{P}(\mu) = \{\gl| D_{cap}(\gl) \; \mbox{matches}
\;D_{wt}(\mu)\}. \ee
Then $\gl \lra f_\gl$ defines a bijection from $\mathbf{P}(\mu)$ to $P(f_\mu).$ Suppose that $\gl \in X^+$  and
      \[\times(\gl) = (c_1, \ldots , c_k)\]
with $c_k < \ldots < c_2 < c_1.$  If $\# \gl = k,$ this means that there are subsets $$\{ i_1< \ldots < i_k\} \subseteq \{1, \ldots
,m\}, \quad \{ j_1> \ldots > j_k\} \subseteq \{1, \ldots ,n\}$$ such
that \[ (\gl  + \gr,  \gep_{i_p}) = (\gl  + \gr,  \gd_{j_p}) =
c_p,\] for $1 \leq  p \leq  k,$ and we set $\ga_p = \gep_{i_p} - \gd_{j_p}.$  Now suppose that the cap in
$D_{cap}(\gl)$ beginning at $c_p$ ends at $d_p = c_p + r'_p.$
Next let
$(r_1,\dots,r_k)$ be the lexicographically smallest tuple of
strictly positive integers such that for all $\theta=
(\theta_1,\dots,\theta_k) \in \{0,1\}^k$, $$\SS_{\theta}(\gl) =  \gl + \sum_{p=1}^k
\theta_p r_p \ga_p \in X_{reg},$$ and   let $\RR_{\theta}(\gl)$ denote the unique element
of $X^+(m|n)$ which is  conjugate under the dot action of $W$ to  $\SS_{\theta}(\gl).$
\bl \label{} We have
\begin{itemize}
\item[{(a)}]
$r_p = r_p'$ for $1 \leq  p  \leq k.$
\item[{(b)}]
$D_{wt}(\RR_{\theta}(\gl))$ is obtained from
$D_{wt}(\gl)$
by interchanging the $\ti$ and $\ci$ located at $c_p$ and  $d_p$ respectively for all  $p$ such that $\theta_p = 1,$ and leaving all other symbols unchanged.
\ei \el \bpf Clearly
\be \label{SSth}  (\SS_{\theta}(\gl)+\gr, \gep_{i_q} ) = (\SS_{\theta}(\gl)+\gr,\gd_{i_q} ) = c_q + \theta_q r_q.\ee
Assume by induction that $r_q = r_q'$ for $1 \leq  q  \leq p-1,$ and set
\[Y_p = \{ c_1, \ldots,c_{p-1}, d_1, \ldots,d_{p-1}\}\cup f_\gl^{-1}(<) \cup f_\gl^{-1}(>).\] From the definition of the cap diagram $D_{cap}(\gl)$ it follows that \[r'_p = \min \{r |r > 0, c_p + r \notin Y_p\}.\]
Using this and (\ref{SSth}) we conclude that $r_p = r_p'$.  This proves (a), and (b) follows since when weights are written as in equation (\ref{alb}), the dot action of $W$ is implemented by permuting the entries.\epf
\bc \label{jbc}  Let $\RR_{\theta}(\gl)$ denote the unique element
of $X^+(m|n)$ which is  conjugate under the dot action of $W$ to $\gl + \sum_{p=1}^k \theta_p r_p \ga_p \in
X_{reg}.$ Then  \begin{itemize}
\item[{(a)}]
$\{\RR_{\theta}(\gl)|\theta \in \{0,1\}^k\} = \{\mu \in
X^+|D_{cap}(\gl)\; \mbox{matches the cap diagram} \;D_{wt}(\mu)\}.$
\item[{(b)}] $\mu = \RR_{\theta}(\gl)$ for some $\theta \in
\{0,1\}^k$ if and only if $\gl \in \mathbf{P}(\mu).$
\end{itemize}
\ec \bpf This  follows at once from the Lemma and equation  (\ref{pmgl}).\epf
The following
reformulation of the main theorem in \cite{Br } was shown to the
first author by Jon Brundan. The result is also
recorded in [Theorem 5.2, BS08a] using slightly different notation.
\bt \label{Jon}  In the Grothendieck group
of the category $\cF$ we have
$$K(\mu) = \sum_{\gl \in \mathbf{P}(\mu)} L(\gl).$$\et \bpf The Main Theorem
in \cite{Br} states that for each $\mu \in X^+(m|n)$,
$$
[K(\mu):L(\gl)] = \left\{
\begin{array}{ll}
1&\hbox{if $\mu = \RR_{\theta}(\gl)$ for some $\theta =
(\theta_1,\dots,\theta_k)
\in \{0,1\}^k$,}\\
0&\hbox{otherwise.}
\end{array}
\right.
$$
The result now follows immediately by Corollary \ref{jbc}. \epf
We now state the main result of \cite{S2} in terms of 
diagrams.
\bt\label{vt} If
$$(1+\sigma_1)\dots (1+\sigma_k)f_{\lambda}=\sum_{\mu} c_{\lambda,\mu}
f_{\mu},$$
then in the Grothendieck group of the category $\cF$ we have
$$K(\lambda)=\sum_{\mu} c_{\lambda,\mu} L(\mu).$$
\et

Combining Theorem A, Theorem \ref{Jon} and Theorem \ref{vt}, we obtain a combinatorial proof of the equivalence of the algorithms from \cite{Br} and \cite{S2}.

In the rest of this section we explain how to deduce Theorem \ref{vt}, and some further results that we will require in Section \ref{seae}, from results in \cite{S2}.
An equivalence of categories allows us to focus our attention on the category $\cF^k$ of all finite dimensional modules which have the degree of atypicality $k$, see \cite{S4}
or \cite{GS}. From now on we will make this assumption.
Let $F^k$ be the set of core-free $f \in F$ such that $\#f = k$.  As before we may identify $\mathbb Z F^k$ with the Grothendieck group of $\cF^k.$\\
\\
 Let $\Delta$ be the set of roots of $\fg$, and for any $\alpha\in\Delta$ denote by $\fg^{\alpha}\subset\fg$ the corresponding root space.
 Let $\gamma\in\fh^*$ be a weight such that
$(\alpha,\gamma)\geq 0$ for all positive roots $\alpha$. Set
$$\Delta_\gamma=\{\alpha\in\Delta | (\alpha,\gamma)\geq 0\}.$$
We say that $\gamma$ {\it defines the
parabolic subalgebra} ${\fq}\subseteq\fg$ where
$${\fq}=\fh\oplus \bigoplus_{\alpha\in \Delta_\gamma}\fg^{\alpha}.$$
 Note that ${\fb}\subseteq\fq$.
Let
$\fl$ be the ad-$\fh$ stable  Levi subalgebra of $\fq$ and note that $\fl$ has a $\mathbb Z$-grading
$\fl=\fl_{-1}\oplus\fl_0\oplus \fl_1$, similar to the $\mathbb Z$-grading of
$\fg$.
Every ${\fl}$-module
can be made into a ${\fq}$-module with trivial action of the
nilpotent radical of ${\fq}.$ In particular,
 the $\mathbb Z$-grading
on $\fl$ allows us to construct a Kac module $K_{{\sfq}}(\lambda)$ for  ${\fq}$ and we denote the unique simple factor module of $K_{\sfq}(\gl)$ by
$L_{{\sfq}}(\lambda)$.
Note that one can write
\be \label{NE2}\ch L_{\sfq}(\gl)=\sum_{\mu\leq\gl}a_{\sfq}(\gl,\mu)\ch K_{\sfq}(\mu), \ee
and we denote the matrix
with coefficients $a_{\sfq}(\gl,\mu)$ by   $A_{\sfq}$.

Now let ${\fq}\subset\fp$ be a pair of parabolic subalgebras, and $V$ be a
finite dimensional ${\sfq}$-module. Let $\Gamma_{\sfp,{\sfq}}(V)$ be the
maximal finite dimensional quotient of the induced module
$U(\fp)\otimes_{U({\sfq})}V$. Then clearly $\Gamma_{\sfp,{\sfq}}$ is a
functor from the category of finite dimensional ${\fq}$-modules to the
category of finite dimensional $\fp$-modules and this functor is exact
on the right. In general,
the functor $\Gamma_{\sfp,\sfq}$ is not exact, but
it was proven in
\cite{S2} that $\Gamma_{\sfp,\sfq}$ is exact on $\fq$-modules that are
free over $U(\fl_{-1})$.
It is not hard to see that any $\fq$-module that is free over
$U(\fl_{-1})$ has a filtration
with quotients isomorphic to Kac modules $K_{\sfq} (\mu)$. (By
definition, $K_{\sfq}(\mu)=K_{\sfl}(\mu)$ with trivial action of the
nilpotent radical of $\fq$). Moreover,
\begin{equation}\label{KM}
\Gamma_{\sfp,\sfq}K_{\sfq}(\mu)=K_{\sfp}(\mu).
\end{equation}

We construct derived functors $\Gamma^i_{\sfp,\sfq}(V)$,
of  $\Gamma_{\sfp,\sfq}$ as follows. Take a resolution
$$\dots\to M^1\to M^0\to 0$$ of $V$ by
$\fq$-modules that are free over $U(\fl_{-1}),$
and define $\Gamma^i_{\sfp,\sfq}(V)$ to be the $i^{th}$ cohomology
group of the complex $$\dots\to \Gamma_{\sfp,\sfq}(M^1)\to  \Gamma_{\sfp,\sfq}(M^0)\to 0.$$
The result does not depend on a choice of resolution since
$\Gamma_{\sfp,\sfq}$ is exact on $\fq$-modules which are free over
$U(\fl_{-1})$.
Clearly,
we have a natural surjective homomorphism of $\fp$-modules
$\gamma:\Gamma^0_{\sfp,\sfq}(L_{\sfq}(\gl))\to L_{\sfp}(\gl)$.
Define
$$U^i_{\sfp,\sfq}(\gl)=\Gamma^{i-1}_{\sfp,\sfq}(L_{\sfq}(\gl))$$
for $i>1$ and
$$U^1_{\sfp,\sfq}(\gl)=\operatorname{Ker}\gamma.$$
Put
$$U^i_{\sfp,\sfq}(\gl,\mu)=[U^i_{\sfp,\sfq}(\gl):L_{\sfp}(\mu)].$$

Take a resolution $M^{\bullet}_\lambda$ of $L_{\sfq}(\gl)$ such that
$M_\lambda^0=K_{\sfq}(\gl)$,
and for $i>0$,
$M_\lambda^i$ is free over $U(\fl_{-1})$
and has all weights strictly less
than $\lambda$.  Then  we have
\begin{equation}\label{lt}
U^i_{\sfp,\sfq}(\gl,\mu)\neq 0\,\, \text{implies}\,\,\gl>\mu.
\end{equation}
Clearly, we have
\be \label{NE3}\ch L_{\sfq}(\gl)=\sum_{i \geq 0}(-1)^i\ch M_\gl^i,\ee
and
\be \label{NE4}\ch L_{\sfp}(\gl)-\sum_{i \geq 1}(-1)^i\ch U^i_{\sfp,\sfq}(\gl)=\sum_{i \geq 0}(-1)^i\ch
\Gamma_{\sfp,\sfq}(M_\gl^i).\ee
Combine  (\ref{NE2}) and (\ref{NE3}), and then apply $\Gamma_{\sfp,\sfq}$, using (\ref{KM}) to obtain
\[\sum_{i \geq 0}(-1)^i\ch
\Gamma_{\sfp,\sfq}(M_\gl^i) = \sum_{\mu\leq\gl}a_{\sfq}(\gl,\mu)\ch K_{\sfp}(\mu).\]
From this and (\ref{NE4}) we deduce the following important identity
\begin{equation}\label{add1}
\ch L_{\sfp}(\gl)-\sum_{\nu,i}(-1)^i U^i_{\sfp,\sfq}(\gl,\nu)\ch L_{\sfp}(\nu)=
\sum_{\mu\leq\gl}a_{\sfq}(\gl,\mu)\ch K_{\sfp}(\mu).
\end{equation}

Set $U_{\sfp,\sfq}(\gl,\mu) = \sum_{i \geq 1}(-1)^i U^i_{\sfp,\sfq}(\gl,\mu)$.
Let $U_{\sfp,\sfq}$ be the matrix
with coefficients  $U_{\sfp,\sfq}(\gl,\mu)$. Then using (\ref{NE2}),
the identity (\ref{add1}) can be rewritten in the form

\begin{equation}\label{add2}
(1-U_{\sfp,\sfq})A_{\sfp}=A_{\sfq}.
\end{equation}
For $1\leq s\leq k$  let
$$\gamma_s=s\varepsilon_1+\dots +\varepsilon_s+s\delta_k+\dots +\delta_{k-s+1},$$
and let ${\fq}^{(s)}$ be the parabolic subalgebra defined by $\gamma_s$.
Consider the flag of parabolic subalgebras \be \label{NE7}\sfg=\fq^{(0)}\supset\fq^{(1)}\supset\dots\supset
\fq^{(k)}=\fb.\ee
Consecutive application of (\ref{add2}) to the pairs
$\fp=\fq^{(i)}\supset\fq=\fq^{(i+1)}$ and the fact that $A_{\sfb}=1$  give us
\begin{equation}\label{add3}
A_{\sfg} =(1-U_{\sfq^{(0)},\sfq^{(1)}})^{-1}\dots(1-U_{\sfq^{(k-1)},\sfq^{(k)}})^{-1}.
\end{equation}
The matrix $C$ with coefficients $c_{\gl,\mu}$ as in Theorem \ref{vt} is the  inverse of
$A_{\sfg}$.
Hence we have
\begin{equation}\label{add3}
C=(1-U_{\sfq^{(k-1)},\sfq^{(k)}})\dots(1-U_{\sfq^{(0)},\sfq^{(1)}}).
\end{equation}

Note that this is an equality of operators on the Grothendieck group of the category $\cF^k$. We define analogous
linear operators
$U_{\sfp,\sfq}$ and $C$ on $\mathbb Z F^k$
by first setting
$$U_{\sfp,\sfq}(f,g)=U_{\sfp,\sfq}(\gl(f),\gl(g)), \quad c(f,g)=c(\gl(f),\gl(g)).$$ and then
\begin{equation}\label{j22}
U_{\sfp,\sfq}(f)=\sum_{g\in F}U_{\sfp,\sfq}(f,g) g,\quad
C(f)=c(f,g)g.
\end{equation}
Then (\ref{add3}) can be also be viewed as an equality of linear operators on $\mathbb Z F^k$.

The equation (\ref{add3}) reduces the problem of finding the composition
factors of Kac modules to the problem of calculating
$U^{i}_{\sfq^{(j)},\sfq^{(j+1)}}({\gl},{\mu})$.
Concerning the latter problem, the next result summarizes Theorems 6.15 and 6.24 from \cite{S2}.
\bt \label{mainselecta}  \begin{itemize}
\item[{}]
 \item[{(a)}] If $\gl - \ga$ is $\fq^{(j)}$-dominant then
\be \label{neweq1}
U^{i}_{\sfq^{(j)},\sfq^{(j+1)}}({\gl}) = U^{i+1}_{\sfq^{(j)},\sfq^{(j+1)}}({\gl-\ga})\nonumber \ee
 for $i > 1,$ and
\be \label{neweq2}
U^{1}_{\sfq^{(j)},\sfq^{(j+1)}}({\gl}) = L_{\sfq^{(j)}}(\gl-\ga) \oplus U^{2}_{\sfq^{(j)},\sfq^{(j+1)}}(\gl-\ga). \nonumber \ee
 \item[{(b)}] If $\gl - \ga$ is not  $\fq^{(j)}$-dominant then \be
   \label{neweq3}  [U^{i}_{\sfq^{(j)},\sfq^{(j+1)}}({\gl}):L_{\sfq^{(j)}}(\mu)] =
   [U^{i-1}_{\sfq^{(j+1)},\sfq^{(j+2)}}({\gl-\ga}):L_{\sfq^{(j+1)}}(\mu)]. \nonumber \ee and
   $U^{1}_{\sfq^{(j)},\sfq^{(j+1)}}({\gl}) = 0.$
\item[{(c)}]\be\label{add5}
U^1_{\sfq^{(k-1)},\sfq^{(k)}}(\gl)=L_{\sfq^{(k-1)}}(\gl-\ga),\quad
U^i_{\sfq^{(k-1)},\sfq^{(k)}}(\gl)=0, \, \text {if} \,\,  i>1.\nonumber \ee
\end{itemize}
\et

To prove Theorem \ref{vt} it remains to interpret the above result in  terms of diagrams.
Below  we use an induction argument  on $k$, and for this purpose we define, by analogy with
(\ref{NE7}), the flag of parabolic subalgebras in $\fgl(k-1,k-1)$
$$\fgl(k-1,k-1)={\fp}^{(0)}\supset {\fp}^{(1)}\dots\supset
{\fp}^{(k-1)}=\fb',$$ where $\fb'$ is the Borel subalgebra consisting of upper triangular matrices. Let $\fg_{(s)}$ be the subalgebra of $\fg$ consisting of all matrices with zero entries in rows $s, 2k-s+1$, and zero entries in columns $s, 2k-s+1$.  We have an obvious isomorphism from $\fgl(k-1,k-1)$ to $\fg_{(s)}$, and we denote the image of ${\fp}^{(j)}$ under this isomorphism by ${\fp}^{(j)}_{(s)}.$
If $l(\fp)$ is the quotient of $\fp$ by the nilpotent radical, then because we deleted two diagonal entries from $\fg$ to get $\fg_{(s)}$ we have
$$l(\fq^{(j+1)})\simeq l(\fp^{(j)}_{(s)})\oplus\mathbb C \oplus\mathbb C.$$

Let $\gl=(a_1,\ldots,a_k|a_k,\ldots,a_1)$ and
$${\gl}'=(a_1,\ldots,a_{s-1},a_{s+1},\ldots,a_k|a_k,\ldots,a_{s+1},a_{s-1},\ldots,a_1).$$
If we regard
$L_{\sfq^{(j+1)}}(\gl)$ as a ${\sfp^{(j)}_{(s)}}$-module, via the above isomorphism, then it remains irreducible with highest weight $\gl'$.
This implies
\be\label{marks}
[U^{i}_{\sfq^{(j+1)},\sfq^{(j+2)}}(\gl):L_{\sfq^{(j+1)}}(\mu)]=
[U^{i}_{\sfp^{(j)},\sfp^{(j+1)}}({\gl}'):L_{\sfp^{(j)}}(\mu')].
\ee

\bl\label{addleg} Let
  $U^{i}_{\sfq^{(j)},\sfq^{(j+1)}}(f)=U^{i}_{\sfq^{(j)},\sfq^{(j+1)}}({\gl}(f))$ and
$L_{\sfq}(f)=L_{\sfq}(\gl(f))$. Next let \be \label{fti}
  f^{-1}(\times)=\{a_1,\dots,a_k\}\ee with $a_1>a_2>\dots > a_k$ and $a=a_{j+1}$.
Then the relations of Theorem \ref{mainselecta} can
  be rewritten in the following way in terms of weight diagrams.
\begin{itemize}
 \item[{(a)}] If $f(a-1)=\ci$  then
\be \label{addeq1}
U^{i}_{\sfq^{(j)},\sfq^{(j+1)}}(f) = U^{i+1}_{\sfq^{(j)},\sfq^{(j+1)}}(f^a_{a-1})\ee
 for $i > 1,$ and
\be \label{addeq2}
U^{1}_{\sfq^{(j)},\sfq^{(j+1)}}(f) = L_{\sfq^{(j)}}(f^a_{a-1})) \oplus U^{2}_{\sfq^{(j)},\sfq^{(j+1)}}(f^a_{a-1}).\ee
 \item[{(b)}] If $f(a-1)=\times$  \be
   \label{addeq3}  [U^{i}_{\sfq^{(j)},\sfq^{(j+1)}}(f):L_{\sfq^{(j)}}(g)] =
   [U^{i-1}_{\sfp^{(j)},\sfp^{(j+1)}}(f^{a}):L_{\sfq^{(j)}}(g^{a-1})]. \ee

In addition,
   $U^{1}_{\sfq^{(j)},\sfq^{(j+1)}}(f) = 0.$
\item[{(c)}]\be\label{add6}
U^1_{\sfq^{(k-1)},\sfq^{(k)}}(f)=L_{\sfq^{(k-1)}}(f^a_{a-1}),\,\;\;
U^i_{\sfq^{(k-1)},\sfq^{(k)}}(f)=0, \, \text {if} \,\,  i>1.\nonumber \ee
\end{itemize}
\el
\begin{proof}
(a) and (c) follow immediately from the identity
$\gl(f^a_{a-1})=\gl(f)-\ga$. To prove (b) note that by  (\ref{marks})
we have
\be\label{NE9}[U^{i-1}_{\sfq^{(j+1)},\sfq^{(j+2)}}(\gl(f)-\ga):L_{\sfq^{(j+1)}}(g)]=
[U^{i-1}_{\sfp^{(j)},\sfp^{(j+1)}}(f^{a}):L_{\sfp^{(j)}}(g^{a-1})].\ee
   By  Theorem \ref{mainselecta} (b) we have
$$[U^{i}_{\sfq^{(j)},\sfq^{(j+1)}}({\gl}(f)):L_{\sfq^{(j)}}(\gl(g))] =
   [U^{i-1}_{\sfq^{(j+1)},\sfq^{(j+2)}}({\gl(f)-\ga}):L_{\sfq^{(j+1)}}(\gl(g))],$$
and combining this with (\ref{NE9}) we deduce (\ref{addeq3}).
\end{proof}
We introduce two related pieces of notation.  Sometimes one is more convenient than the other. First suppose $f \in F^k$ with $f(a) = \ti,$ set $$ \lm_k(f,a,i) = \{b \in  \mathbb{Z}| f^a_b \; \mbox{ is obtained from} \;  f \; \mbox{by a legal move of weight} \;i \}. $$
Next define with the notation of (\ref{fti}),
$$\LM(f,p) = \{g|g \; \mbox{ is obtained from} \;  f \; \mbox{by a legal move of weight} \; 0
\;  f \; \mbox{starting at} \;a_p
\}.$$
\bl \label{OMT}  Let $f\in F^k$, $f(a)=\ti$, then we have
\begin{itemize}
\item[{(a)}] If $f(a-1) = \ci$, and $i > 0,$ then  $\lm_k(f,a,i) = \lm_k(f^a_{a-1},a-1,i+1).$
\item[{(b)}] If $f(a-1) = \ci$,  then  $\lm_k(f,a,0) = \lm_k(f^a_{a-1},a-1,1)\cup \{{a-1}\}.$
 \item[{(c)}] If $f(a-1) = \ti$, and $h =f^a,$ then  $\lm_k(f,a,0) = \emptyset$ and $$\lm_k(f,a,i) = \lm_{k-1}(h,a-1,i-1),$$ for $i>0$.
 \item[{(d)}] If $f(a-1) = \ci$, and $h =f^a_{a-1},$ or $f(a-1) = \ti$, and $h =f^a,$ then $f^a_b = h_b^{a-1}$
 for all $b \in \lm_k(f,a,i).$

\end{itemize}
\el \bpf Straightforward. \epf
\bc\label{addleg1} Let
$U^{i+1}_{\sfq^{(j)},\sfq^{(j+1)}}(f)=\sum_{g\in F}
U^{i+1}_{\sfq^{(j)},\sfq^{(j+1)}}(\gl(f),\gl(g))g$.
Then
\be U^{i+1}_{\sfq^{(j)},\sfq^{(j+1)}}(f) = \sum_{b \in \lm_k(f,a,i)}f^a_b.
\nonumber \ee \ec
\begin{proof}
It is sufficient to prove the statement
for $j=0,$ since none of the terms in $U^{i+1}_{\sfq^{(j)},\sfq^{(j+1)}}(f)$ depends on the $j$ rightmost $\ti$-s in $D_{wt}(f)$.\\
\\
The proof goes by induction on $k$, and for $k$ fixed a second induction on the distance between the leftmost $\times$ and the
rightmost $\times$ of $D_{wt}(f)$.
The case $k=1$ immediately follows from Lemma \ref{addleg}(c).
Let $a=a_1$ be the position of the rightmost $\times$ of $D_{wt}(f)$.

First, assume that $f(a-1)=\ci$. Let $h=f^a_{a-1}$. Then  if $i > 0$, we have using parts (a), (d) of  Lemma \ref{OMT}, the second induction hypothesis applied to $h$, and then (\ref{addeq1}) we have
\begin{eqnarray}\label{}
\sum_{b \in \lm_k(f,a,i)}f^a_b &=&\sum_{b \in \lm_k(h,a-1,i+1)}
h_b^{a-1}\nonumber \\
&=&\sum_g
U^{i+2}_{\sfq^{(j)},\sfq^{(j+1)}}(h, g)g \nonumber \\ &=& \sum_{g}
U^{i+1}_{\sfq^{(j)},\sfq^{(j+1)}}(f, g)g. \nonumber
\end{eqnarray}  If $i = 0,$ the result follows similarly, using part (b) of the Lemma and  (\ref{addeq2}).  Finally if $f(a-1)=\times$, let $h=f^a$. Then using parts (c) , (d) of the Lemma, induction on $k$ and (\ref{addeq3}) we have
\begin{eqnarray}\label{}
\sum_{b \in \lm_k(f,a,i)}f^a_b &=&\sum_{b \in \lm_{k-1}(h,a-1,i-1)}
h_b^{a-1} \nonumber \\ &=& U^{i}_{\sfp^{(j)},\sfp^{(j+1)}}(h) =
U^{i+1}_{\sfq^{(j)},\sfq^{(j+1)}}(f). \nonumber
\end{eqnarray}
\end{proof}

\bc\label{sigma}
We have
\be \label{NE8}\sigma_{j+1}=-U_{\sfq^{(j)},\sfq^{(j+1)}}.\ee
\ec
\begin{proof} The result follows from the previous Corollary since
\[U_{\sfq^{(j)},\sfq^{(j+1)}}= \sum_{i}(-1)^i U_{\sfq^{(j)},\sfq^{(j+1)}}.\]
\end{proof}
Together Equations (\ref{add3}) and (\ref{NE8}) yield Theorem \ref{vt}.

We remark that by Corollary 6.25 from \cite{S2},
the modules $U^{i}_{\sfq^{(j)},{\sfq^{(j+1)}}}({\gl})$ are semisimple.  Thus Corollary \ref{addleg1} determines their decompositions into simple modules.  In particular this gives us the first part of the next result. The second part will be used in Section 6 of this paper.

\bc\label{selecta2}
For  $f$ be as in {\rm (\ref{e2})} we have
\begin{itemize}
\item[{(a)}] \be U^{1}_{\sfq^{(p-1)},
{\sfq^{(p)}}}(\lambda(f)) =
\bigoplus_{g\in \LM(f,p)}L_{{\sfq}^{(p-1)}}(\lambda(g)). \nonumber \ee
\item[{(b)}]
$\Gamma_{\sfq^{(p-1)},\sfq^{(p)}}(L_{{\sfq}^{(p)}}(\lambda(f)))$ is
generated by a highest weight vector of weight $\lambda$ and its
structure can be described by the exact sequence
\begin{equation}\label{sec}
0\to \bigoplus_{g\in \LM(f,p)}L_{{\sfq}^{(p-1)}}(\lambda(g)) \to
\Gamma_{\sfq^{(p-1)},\sfq^{(p)}}(L_{{\sfq}^{(p)}}(\lambda(f)))\to
L_{{\sfq}^{(p-1)}}(\lambda(f))\to 0.\nonumber
\end{equation}\ei
\ec \bpf  By what we said above, it is enough to note that (b) follows from (a) and Lemma 4.11 in \cite{S2}.\epf
 \section{The Graph $\mathcal{G}$ and the Involution on Irregular Paths.}
From now on we consider only elements of $F$ that are core-free.
Define  $\mathcal{G}$ to be the  oriented graph whose vertices are elements of
$F$, and we join $f$ and $g$ by an edge $f\longrightarrow g$ if $g$
is obtained from $f$ by a  legal move.   We put the label $[s,t]$ on
this edge if the corresponding legal move has start $s$ and end $t$,
in other words, $g=f_t^s$ (always $s>t$). The weight of an edge is the
weight of the corresponding legal move, and if $g=f_t^s$ as above we set $l([t,s]) = l_f(t,s)$ \\
\\
  It is easy to check that
$\mathcal{G}$ does not have oriented loops.
A {\it path}  in $\mathcal{G}$  is a sequence
$[s_1,t_1],\dots,[s_q,t_q]$ where for $1 \leq  i \leq q$,
$[s_i,t_i]$ is a legal move from $f_{i-1}$ to $f_i$.
We say that the path is
 {\it increasing} if $s_1<\ldots<s_q$. (It follows immediately from the definition that
in any path $s_i\neq s_{i+1}$.) Often we refer to a path by listing only the legal moves. The weight $l(P)$ of a path $P$ is
the sum of weights of all edges in $P$.

\bl\label{weak} Let $\cP_{f,g}(\mathcal{G})$ denote the set of all
increasing paths in $\mathcal{G}$ leading from $f$ to $g$, and let
$$(1 + \sigma_1) \ldots (1 + \sigma_k)f =\sum_{g} c_{f,g} g.$$
Then \be\label{first}
c_{f,g}=\sum_{P\in\cP_{f,g}(\mathcal{G})}(-1)^{l(P)}. \ee \el
\begin{proof} Write
$$(1 + \sigma_1) \ldots (1 + \sigma_k)f = \sum_{i_1<\ldots<i_r}\sigma_{i_1}\ldots\sigma_{i_r}(f).$$
Using (\ref{sif}) we see that each increasing path $P$ with edges
$[a_{i_r},b_{i_r}],\ldots,[a_{i_1},b_{i_1}]$ which leads from $f$ to
$g$ gives the term  $(-1)^{l(P)}g$ in
$\sigma_{i_1}\ldots\sigma_{i_r}(f)$.
\end{proof}

We call an increasing path from $f$ to $g$ in $\mathcal{G}$
{\it  irregular} if one of the following conditions hold
\begin{itemize}
\item[{(a)}] The path contains an edge with positive weight
\item[{(b)}]
  There are repetitions among the labels on the path, in other words
there are edges with label $[c,d]$ and $[b,c]$ in the path. \ei
 An edge $[c,d]$
of an irregular path is called {\it  irregular} if it has a positive
weight or there is an edge with label $[b,c]$ later in the path. A
path which is not irregular is {\it  regular}.

\bl \label{medium} Suppose an increasing path  has edges $[b,c]$ and
$[a,s]$ with $c < s < b < a.$  Then the  edge $[b,c]$ has positive
weight. \el \bpf If the result is false, then with $b, c$ fixed
choose a counterexample with $s$ as large as possible.
Suppose that
the edge $[b,c]$ connects vertex $f$ to $f'$ and  that $D =
D_{wt}(f')$ has $P \; \ti'$s and $p \; \ci'$s in the interval
$(c,s)$. Similarly suppose that the edge $[a,s]$ connects vertex
$g'$ to $g$ and  that $D_{wt}(g')$ has $Q \; \ti'$s and $q \; \ci'$s
in the interval $(s,b)$.
We claim that $D$ also has $Q \; \ti'$s and $q \; \ci'$s in the interval $(s,b)$.  Indeed, consider the part of the path between the edges $[b,c]$ and
$[a,s].$  Since the path is increasing no $\ti$ in the interval $(s,b)$ can be moved.  Also by the choice of the counterexample, there can be no edge with label $[d,e]$ where $d > b$ and $s<e<b.$ The claim follows from this. We deduce that $P \geq p$ and $Q \geq q,$ since $[b,c]$ and $[a,s]$ are legal moves.  Now $D$ has $p+q+1 \; \ci'$s in the interval $(c,b)$ since, in addition to those counted before there is also a $\ci$ at $s$.  Because $[b,c]$ has weight zero,  and $D$ has $P+Q \; \ti'$s in the interval $(c,b)$, we have $P+Q = p +q+1$.  This implies that either $P=p,$ in which case the cap in $D_{cap}(f')$ beginning at $c$ would end at $s$, or $Q=q,$ in which case the cap in $D_{cap}(g')$ beginning at $s$ would end at $b$.  Either way we reach a contradiction.\epf

\bl\label{strong} Let
$\cR_{f,g}({\mathcal{G}})$ denote the set of all increasing regular
paths in ${\mathcal{G}}$ leading from $f$ to $g$. Then \be\label{second}
\sum_{P\in\cP_{f,g}(\mathcal{G})}(-1)^{l(P)}=|\cR_{f,g}({\mathcal{G}})|.
\ee \el
\begin{proof} 
Define an involution $*$ on the set of all increasing irregular
paths by the following procedure. Let $P$ be an irregular increasing
path. Consider  the irregular edge $[s,t]$ of $P$ with maximal
possible end $t$. 
There are exactly two possibilities: either $P$ contains a
regular edge $[a,s]$, or $s$ is not the end of any edge in $P$.\\ \\
In the former case, define $P^*$ to be the path obtained from $P$ by
removing $[s,t]$ and $[a,s]$ and inserting the edge $[a,t]$.
If there were an edge $[b,c]$ in $P$ with $s<b<a$ and $t<c<s$, then  $[b,c]$ would be irregular by Lemma \ref{medium}, contradicting the choice of $s$.
Since the edge $[a,s]$ is regular, it has zero weight.
Therefore $l([a,t])=l([s,t])+1$ and
$l(P^*)=l(P)+1$. Note also that $P^*$ is again irregular and the
edge $[a,t]$ is the irregular edge with maximal possible end.\\
\\
In the latter case, let  $f'\to g'$ be the  edge with label $[s,t]$. Then
$g'(s)=\ci$, and the symbol $\ti$ occurs more often than $\ci$ in
the part of the weight diagram $D_{wt}(g')$ strictly between $t$ and
$s$. In other words we can find a cap in $D_{cap}(g')$ beginning at
$b>t$ and ending at $s$. Then we define $P^*$ to be the path
obtained from $P$ by removing the edge $[s,t]$ and inserting the
edges $[b,t]$ and $[s,b]$. Note that $[s,b]$ is regular, $[b,t]$
is irregular and $l([b,t])=l([s,t])-1$. Hence $l(P^*)=l(P)-1$. It is
clear that $P^*$ is
irregular and $[b,t]$ is the irregular edge with
maximal possible end.

It is obvious that $*$ is an involution and since
$(-1)^{l(P)}=-(-1)^{l(P^*) }$,  all irregular paths in the left hand
side of (\ref{second}) cancel. Hence we have
$$\sum_{P\in\cP_{f,g}(\mathcal{G})}(-1)^{l(P)}=\sum_{P\in\cR_{f,g}({\mathcal{G}})}(-1)^{l(P)}.$$
Now the statement follows since $l(P)=0$ for any regular
path $P$.

\end{proof} \noi We have two immediate consequences of the above work, namely
\be\label{first1} c_{f,g}=|\cR_{f,g}({\mathcal{G}})|, \ee and
\be\label{strong1}(1 + \sigma_1) \ldots (1 + \sigma_k)f =\sum_{g}
|\cR_{f,g}({\mathcal{G}})|g. \ee \bl \label{} Suppose $f, g \in F.$
\begin{itemize}
\item[{(a)}] If $g \in P(f),$ then  $|\cR_{f,g}({\mathcal{G}})| = 1.$ \item[{(b)}] If $g \notin P(f)$
then $\cR_{f,g}({\mathcal{G}})$ is empty.
\end{itemize}
\el \bpf Suppose $g \in P(f),$ and that $\ti(f) = (a_1, a_2,
\ldots,a_k).$  Then let $$I = \{ i \in \{1,\ldots,k\}|D_{cap}(g)\;
\mbox{has a cap ending at} \; a_i\},$$  and for $i \in I,$ suppose
that ending at $a_i$ begins at $b_i$. Then
there is a regular path from $f$ to $g$ given by
${\stackrel{\longrightarrow}{{\rm \prod}}}_{i\in I}[a_i,b_i]$ where
the arrow means that we take the product in the order that gives an
increasing path. It follows easily from the definitions that this is the only way to get a regular
increasing path from $f$ to $g$. \epf
Theorem A immediately follows from (\ref{strong1}) and  the previous Lemma.
\bexa {\rm  Let $k = 2.$ For $a < b \in \mathbb{Z}$ define
$f_{(a,b)} \in F$ so that $\ti(f_{(a,b)}) = (a,b).$ Below we give
the part of the graph $\mathcal{G}$ used to show that \be \label{examp} (1 +
\sigma_1)(1 + \sigma_2)f_{(2,3)} = f_{(2,3)} + f_{(1,3)} +
f_{(0,1)}.\ee Legal moves are represented by arrows together with
their labels.  All edges have weight zero except the edge with label
$[3,1]$ which has weight 1.

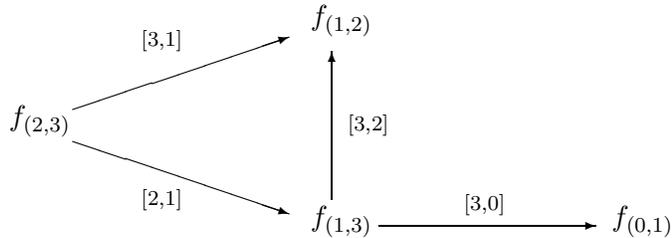
\begin{picture}(258,126)(-55,10)
\thinlines \put(14.0,82.0){$f_{(2,3)}$}
\put(128,120){$f_{(1,2)}$}
\put(128,44){$f_{(1,3)}$}
\put(242,44){$f_{(0,1)}$}
\put(38.0,76.0){\vector(3,-1){82}}
\put(38.0,88.0){\vector(3,1){82}}
\put(154,44){\vector(1,0){82}}
\put(136.0,54.0){\vector(0,1){56}}
\put(64,112){{$\scriptstyle [3,1]$}}
\put(64,52){$\scriptstyle [2,1]$}
\put(142,80){$\scriptstyle [3,2]$}
\put(186,50){$\scriptstyle [3,0]$}
\end{picture}

\noi  There are two irregular paths starting from $f_{(2,3)}$, both
ending at $f_{(1,2)}$. These paths are interchanged by the
involution *.  Summing over the remaining paths and using
(\ref{strong1}) gives (\ref{examp})}. \eexa

\section{Composition factors of Kac modules.}
We describe a procedure for determining the composition factors of
any Kac module, without cancelation.
By Brundan's Theorem we need a procedure for finding the set in
$P(f)$ in Equation (\ref{Pfdef}). This is a problem in
enumerative combinatorics, similar to the problem of describing the
set $B_n$ in Example \ref{catalan}. We give a solution based on the
Lemma below.
By Remark \ref{rem2.1} we can restrict our attention to the core-free
case.

Suppose $f \in F$ with $\ti(f)$ as in Equation
(\ref{e2}) and set $a = a_1$.  Let $f' \in F$ be given by $\ti(f') =(a_2, \ldots, a_{k})$, and set
\[ P(f') = \{\;g' \in F \;| \; D_{cap}(g') \; \mbox{matches} \; D_{wt}(f')  \}, \]
\[ Q(f) = \{\;g \in P(f) \;| \; D_{cap}(g) \; \mbox{has a cap joining} \;  a
\; \mbox{to} \; a +1 \}. \] \bl \label{Lemma 1} There is a bijection
$P(f') \longrightarrow Q(f)$ such that $g' \in P (f')$
maps to $g$ where $\ti(g) = (a,b_1, \ldots, b_{k-1})$ if $\ti(g') =
(b_1, \ldots , b_{k-1})$. \el \bpf Straightforward. \epf Given $f$
as above, we can assume by induction that we have found
$P(f')$ and hence $Q(f)$.  Now suppose $g \in P(f) \backslash
Q(f)$.  Then $D_{cap}(g)$ has a cap joining $b$ to $a$ for some $b <
a$. Replacing this cap with a cap beginning at $a$, we obtain a cap
diagram $D_{cap}(h)$ for some $h \in Q(f)$ such that $g = h^a_b$.
Moreover we can determine the set $P(f) \backslash Q(f)$ as follows.
For each $h \in Q(f)$ list the cap diagrams $D_{cap}(h^a_b)$ that
match $D_{wt}(f)$. Then by
Proposition \ref{mainprop} below, every cap diagram
$D_{cap}(g)$ with $g \in P(f) \backslash Q(f)$ will have been listed
exactly once.\\
\\
The above procedure suggests another  proof of Theorem A. By
induction we may assume that
\[ (1 + \sigma_1) \ldots (1+\sigma_{k-1})f' = \sum_{h \in P(f')} h . \]
Since $\sigma_2, \ldots, \sigma_{k}$ do not move the rightmost
$\ti$ in $D_{wt}(f)$, it follows that
\[ (1 + \sigma_{2}) \ldots (1 + \sigma_k)f = \sum_{h \in Q(f)} h . \]
Now set $\sigma = \sigma_1$.  It remains to show that
\be \label{ss}   \sigma \sum_{h \in Q(f)} h = \sum_{g \in P(f) \backslash Q(f)} g , \ee
but this follows from the Proposition below.  Note that $g
\in P(f) \backslash Q(f)$ implies that \be \label{el00l} \ti (g) =
(b_1, \ldots, b_k) \ee with $b_1 < a$. For $m \in \mathbb{Z}F$, write $m = \sum_{f \in F} |m:f|f$,
with $|m:f| \in \mathbb{Z}$.
Next suppose  $f$ satisfies (\ref{e2}), and that  $f(a) = \ti,$ and
$f(b)  =  \ci.$
\label{revised}\bp \label{mainprop} Suppose $g$ satisfies Equation (\ref{el00l}),
and $b_1 < a.$ Set
\[R_{f,g} = \{h \in Q(f)| 
\; | \sigma h: g| \neq 0\}. \]  Then one of
the following holds\\
\\
(a) $g \in P(f) \backslash
Q(f)$.  In this case $R_{f,g} = \{f\}$ is a singleton and $|\sigma f: g| = 1$.  \\
\\
(b) $g \not\in P(f)$.  In this case either $R_{f,g}$ is empty or
$R_{f,g} = \{ h,h'\}$ consists of two elements and
\be \label{shg} | \sigma h:g| + |\sigma h' : g| = 0. \ee
\ep \noi Theorem A follows immediately from the Proposition and what we have said above.\\ \\ Note that $| \sigma h: g| \neq 0$ if and only if $g$ is
obtained from $h$ by a legal move.  The strategy to prove the Proposition is to identify some caps that are unaffected by the legal move, and remove these, thus reducing to a special case.   The special case is described in the Lemma below, and we omit the proof since it is easily verified.
\bl \label{spcase} Suppose $\ti(g) = (r-1, \ldots , 0),$ and for $0 \le b \le r-1$ set $h^{(b)} = g^r_b.$  
 \bi \item[{(a)}]  
\quad  If $b =  0$, and $f = h^{(b)},$ then  $R_{f,g} = \{f\}.$
\item[{(b)}]  $\quad $ Suppose that  $0 < b \le r-1$, and $f = h^{(b)}$.  If $h = f$ and $h' = h^{(b-1)}$, then  $R_{f,g} = \{h,h'\}$  and Equation (\ref{shg}) holds.
\item[{(c)}]   \quad If $f \neq h^{(b)}$  for any   $b \in [0,r-1],$ then  $R_{f,g}$ is empty.\ei \el 

\noi To reduce to the special case we need to consider separately the caps in $D_{cap}(g)$ that end before $a$ and those that do not. 
To deal with the former, 
set
$${\bf X} = \{(b(C),e(C))| C \mbox{ is a cap in } D_{cap}(g) \mbox{ with  } e(C) < a\},$$  
\noi and suppose that  ${\bf X} =\{(b_i,e_i)|1 \le i \le s\}$.
This allows us to handle the caps in $D_{cap}(g)$ that end before $a$. For the remaining caps, suppose
$p \in F$ is such that $p(i) = \ci$ for $i > a,$ and set
  \be \label{Yp} {\bf Y}(p) = \{b\in \mathbb{Z} | p(b) = \ti, \mbox{ and  the cap in } D_{cap}(p) \mbox{ with  } b = b(C) \mbox{ has } e(C) \ge a\}.\ee
Then provided $p(a) = \ci$, \[ {\bf Y}(p) = \{b \in \mathbb{X} |p \; \mbox{is
obtained from}\; p_a^b\; \mbox{by a  legal move}\} .\]
 Also we have  a disjoint union
\be \label{disun} \ti(g) = {\bf Y}(g) \cup \{b_i\}_{i=1}^s.\ee
Suppose that ${\bf Y}(g) = \{c_0  < c_2 < \ldots < c_{r-1}\}$, and set $c_r = a$.
\bl \label{r1} If $h \in R_{f,g},$ then 
\be \label{r2} h(b_i) = f(b_i) \mbox{ and } h(e_i) = f(e_i) \mbox{ for  } 1 \le i \le s. \ee
\el
\bpf This follows since $h \in Q(f),$ and the caps in $D_{cap}(g)$ that end before $a$  are also caps in $D_{cap}(h)$. (We note that $D_{cap}(h)$ may have an extra cap joining $c_i$ to $b$ for some $i$.) \epf
\noi {\it Proof of Proposition \ref{mainprop}.} 
 Consider $p \in F$ such that Equation (\ref{r2}) holds when $h$ is replaced by $p$, ${\bf Y}(p) \subset \{c_0  < c_2 < \ldots < c_{r}\}$ and  $p(c) = \ci$ for all $c >a.$
Define   $\bar p \in F$ by $\bar p(i)  = p(c_i)$ for $0 \le i \le r$,  and $\bar p(c) = \ci$, if $c \neq c_i$ for any $i \in [0,r]$. 
We claim there is a bijection $R_{f,g} \lra R_{\bar f, \bar g} $ given by $h \lra \bar h$,  and that 
\be \label{r3} |\sigma h: g| = |\sigma \bar h: \bar g|\ee for $h \in R_{f,g}.$ Indeed if $h,h'\in R_{f,g}$ and  $\bar  h = \bar h',$ then $h(b_i) = h'(b_i)$ for all $i$  by Lemma \ref{r1}, and $h(c_i) = h'(c_i)$ by definition of the map $h \lra \bar h$. This shows our map is injective, and surjectivity is shown similarly.  

Finally by (\ref{Yp}) $|\sigma h: g| \neq 0,$ if and only if $g = h^a_c$ for some $c = c_i \in {\bf Y}(g)$.  Moreover in this case we have $\bar g = \bar h^r_i$, and  both sides of  equation (\ref{r3}) are equal to $r-i-1.$ This proves the claim, and the Proposition now follows from Lemma \ref{spcase}.
\hfill $\Box$\\ \\
 We explain how the
above proof is related to our first proof of Theorem A.
To do this we use  the following result.
\bl \label{unique}  \begin{itemize}
\item[] \item[{(a)}] Given $g \in F$ and $a \in \mathbb{Z}$ there is at most
one $f \in F$ such that there is legal move with weight zero from $f$ to $g$ starting
at $a$. \item[{(b)}] For each $g \in P(f)$ there is a unique regular
path from $f$ to $g$. \ei \el \bpf
(a) If there is a cap $C$  in $D_{cap}(g)$ with
$e(C) = a,$ then the unique $f$ in the statement  is  $g_a^b$ where  $b(C) = b$.\\
(b) follows at once from (a).
\epf
Let
* be the involution of the set of
all increasing paths in ${\mathcal{G}}$ leading from $f$ to $g$ defined
in the proof of Lemma \ref{strong}. Then * preserves the set $S$ of
paths with last label of the form $[a,t]$ all of whose edges are
regular except possibly the last. Suppose the set $R_{f,g}$ is
defined as in the proof of Proposition \ref{mainprop} is nonempty.
If $g \in P(f) \backslash Q(f)$ and $R_{f,g} = \{h\},$ then on the
unique regular path from $f$ to $g$, $h$ is the vertex before $g$
and $g$ is obtained from $h$ by a legal move with start $a$. On the
other hand if $g \notin P(f) ,$ then $R_{f,g} = \{ h,h'\}$
and  $S$ consists of two paths which are interchanged by *. The
vertices before the last in these paths are $h$ and $h'$.

\section{The Subgraph $\mathcal{E}$ and Extensions.}\label{seae}

Let $\mathcal{E}$ be the subgraph of  $\mathcal{G}$ obtained by
deleting all edges of positive weight.
Recall that edges of  $\mathcal{E}$ which are in bijection with legal
moves $f\to g$ of
weight $0$. These are especially easy to understand, since they correspond to
switching the labels at the ends of a cap in $D_{cap}(g)$,
The goal if this section is to prove\\ \\
{\bf Theorem B.} \label{ext}We have \begin{itemize}
\item[{(a)}]
$$\operatorname{dim}\operatorname{Ext}^1(L(\lambda(f)),L(\lambda(g)))\leq  1.$$
\item[{(b)}]
$\operatorname{Ext}^1(L(\lambda(f)),L(\lambda(g)))\neq 0$ if and only
if $f\longrightarrow g$ or $g\longrightarrow f$ is an edge of $\mathcal{E}$.\ei
Let $\tau$ be the automorphism of $\fg$ defined by
$\tau(X)=-X^{st}$, where $X^{st}$ is the supertranspose of $X,$ and for any $\fg$-module $M,$ let $M^{\tau}$ denote
the twist by $\tau$.
Thus as a set $M^{\tau} = \{m^\tau|m \in M\}$ and the module structure is given by \[xm^\tau = (\tau(x)m)^\tau\] for $x \in \fg$ and $m \in M.$ The superHopf algebra structure of $U(\fg)$ allows us to make the dual $N^*$ of any $\mathbb{Z}_2$-graded module $N$ in $\cF$ into a module in $\cF$, and we set $\check{M}=(M^*)^{\tau}$. Then $M\to \check{M}$ is a contravariant exact functor on $\cF$ which maps a simple
finite dimensional module to itself. Hence we have
\begin{equation}\label{duality}
\operatorname{Ext}^1(L(\lambda),L(\mu))=\operatorname{Ext}^1(\check{L}(\mu),\check{L}(\lambda))=\operatorname{Ext}^1(L(\mu),L(\lambda)).
\end{equation}
\br {\rm Equation (\ref{duality}) reflects a more general phenomenon.  Indeed there is  well developed theory of links between prime ideals in a Noetherian ring $R$, see for example, \cite{GW}. The {\it graph of links} is the directed graph whose vertices are the prime ideals of $R$, with  arrows between linked prime ideals.  It is shown in \cite{M} that if $\fg$ is a classical simple Lie superalgebra and $\fg \neq P(n),$ then for prime ideals $P, Q$ of $U(\fg)$ there is a link from $P$ to $Q$ if and only if there is a link from $Q$ to $P$. Equation (\ref{duality}) follows from this fact by taking $P$ and $Q$ to be coartinian. In the case where $\fg = \fsl(2,1)$ graph of links between primitive ideals is described in \cite{M}.} \er

Define an order on the set $X(m,n)$ by putting $\mu\leq\lambda$ if
$\lambda-\mu$ is a sum of positive roots.

\bl\label{ext1} Let $\operatorname{Ext}^1(L(\lambda),L(\mu))\neq 0$,
then either $\lambda\leq\mu$ or $\mu\leq\lambda$.
\el
\begin{proof} Assume  that $\lambda$ and $\mu$ are not compatible. Consider an exact sequence
$$0\to L(\lambda) \to M \to L(\mu) \to 0.$$
Since  $\mu$ has multiplicity one as a weight of $M$, a non-zero vector of weight
$\mu$ generates a proper submodule in $M$. Hence the exact sequence splits.
\end{proof}

\bl\label{par} Let ${\fp}\subset\fq$ be a pair of parabolic subalgebras,
and suppose $\gamma\in\fh^*$ defines ${\fp}$. If $\mu\leq\lambda$ and
$\operatorname{Ext}^1_{\sfq}(L_{\sfq}(\lambda),L_{\sfq}(\mu))\neq 0$, then one of the following holds
\begin{itemize}
\item[{(a)}]
$(\mu,\gamma)=(\lambda,\gamma)$ and
$\operatorname{Ext}^1_{{\sfp}}(L_{{\sfp}}(\lambda),L_{\sfp}(\mu))\neq 0$
\item[{(b)}]
$(\mu,\gamma)<(\lambda,\gamma)$ and $L_{\sfq}(\mu)$ is a subquotient in $\Gamma_{{\sfq},{\sfp}}(L_{\sfp}(\lambda))$.\ei
\el
\begin{proof} The condition  $\mu\leq\lambda$ implies that $(\mu,\gc)\leq (\lambda,\gc)$.
Consider a non-split exact sequence
$$0\to L_{\sfq}(\mu) \to M \to L_{\sfq}(\lambda) \to 0.$$
Let $z\in \fh$ be
the element such that $\beta(z)=(\beta,\gamma)$ for every $\beta\in\fh^*$.
For every ${\fq}$-module $N$ let
$$N'=\{x\in N | z x=(\lambda,\gamma)x\}.$$

If  $(\mu,\gamma)=(\lambda,\gamma)$ we have $L_{\sfq}(\lambda)'=L_{\sfp}(\lambda)$,
$L_{\sfq}(\mu)'=L_{\sfp}(\mu)$ and an exact sequence of ${\fp}$-modules
$$0\to L_{\sfp}(\mu) \to M' \to L_{\sfp}(\lambda) \to 0.$$
We claim that this exact sequence does not split. Indeed, $M$ is a
quotient of $U({\fq})\otimes_{U({\sfp})}M'$. If $M'=L_{\sfp}(\lambda) \oplus L_{\sfp}(\mu)$, then
we have an exact sequence
$$U({\fq})\otimes_{U({\sfp})}L_{\sfp}(\lambda) \oplus
U({\fq})\otimes_{U({\sfp})}L_{\sfp}(\mu)\to M\to 0,$$
hence
$M=L_{\sfq}(\lambda) \oplus L_{\sfq}(\mu)$, a contradiction.

If $(\mu,\gamma)<(\lambda,\gamma)$, then $L_{\sfq}(\lambda)'=L_{\sfp}(\lambda)$ and
$L_{\sfq}(\mu)'=0$. Therefore $M'=L_{\sfp}(\lambda)$. The homomorphism
$M'\to M$ of $\sfp$-modules induces a homomorphism
$U({\fq})\otimes_{U({\sfp})}M'\to M$, which is surjective because
$U({\fq})M'=M$. Thus, $M$ is a quotient of
$U({\fq})\otimes_{U({\sfp})}L_{\sfp}(\lambda)$,
hence of  $\Gamma_{{\sfq},{\sfp}}(L_{\sfp}(\lambda))$.
\end{proof}

Lemma \ref{par}, (\ref{duality}) and Corollary \ref{selecta2} imply the following two corollaries.
\bc\label{corext1} If
$\operatorname{Ext}^1(L(\lambda(f)),L(\lambda(g)))\neq 0$ and
$\lambda(f)<\lambda(g)$,
then there is a legal move of weight zero from $f$ to $g$.
\ec
\bpf We have $\operatorname{Ext}^1(L_{{\sfq}^{(k)}}(\lambda(f)),L_{{\sfq}^{(k)}}(\lambda(g))) = 0$. So if $p$ is chosen minimal with $$\operatorname{Ext}^1(L_{{\sfq}^{(p)}}(\lambda(f)),L_{{\sfq}^{(p)}}(\lambda(g))) = 0,$$ then by Lemma \ref{par} with  ${\sfq} = {\fq}^{(p-1)}$ and ${\fp} ={\fq}^{(p)}$,  $L_{\sfq}(\lambda(g))$ is a subquotient of $\Gamma_{{\sfq},{\sfp}}(L_{\sfp}(\lambda(f))$.  Hence the result follows from Corollary \ref{selecta2}. \epf

\bc\label{corext2}
$\operatorname{dim}\operatorname{Ext}^1(L(\lambda(f)),L(\lambda(g)))\leq 1$.
\ec \bpf There is at most one legal move joining $f$ and $g$. Indeed,
if $g$ is obtained from $f$ by a legal move,then $g=f^a_b$, and we have
$$g(a)=\circ,f(a)=\times, g(b)=\times, f(b)=\circ,$$
and $f(s)=g(s)$ if $s\neq a,b$. In other words, $f$ and $g$ are
different exactly in two positions which define the start and the end
of a legal move.
\epf

\bl\label{ext3} Let $g$ be obtained from $f$ by a legal move of weight
zero with start $s_p$. Then
$$\operatorname{Ext}^1_{{\sfq}^{(p-1)}}(L_{{\sfq}^{(p-1)}}(\lambda(g)),L_{{\sfq}^{(p-1)}}(\lambda(f)))\neq 0.$$
\el
\begin{proof} To simplify notation we set
$$\Gamma^{(p)}=\Gamma_{\sfq^{(p-1)},\sfq^{(p)}}.$$
To construct a non-trivial extension consider the exact
sequence from Lemma \ref{selecta2} (b), and set
$$M:=\frac{\Gamma^{(p)}(L_{{\sfq}^{(p)}}(\lambda(f)))}{
\bigoplus_{\stackrel{h\in  \LM(f,p),} {h\neq g}}L_{{\sfq}^{(p-1)}}(\lambda(h))}.$$ Then $M$
is indecomposable and can be included as the middle term in the exact sequence
$$0\to L_{{\sfq}^{(p-1)}}(\lambda(g))\to M \to
L_{{\sfq}^{(p-1)}}(\lambda(f))\to 0.$$
\end{proof}
\noi For $f$ as in (\ref{e2}), define $|f|= a_1+ \dots +a_k.$
\bl\label{extcomb} Let $g\in \LM(f,p), h\in \LM(f,r)$ and $p<r$. \begin{itemize}
\item[{(a)}]  $|f| - |g| \equiv 1$ {\rm mod 2.}
\item[{(b)}]
$g$
and $h$ are not connected by a legal move of weight zero.\ei
\el
\begin{proof} (a) follows immediately from the definition of a legal move.  For (b), assume the opposite. Since $h(s_p)=\times$ and
  $g(s_p)=\circ$ the only possibility is $g\in \LM(h,p)$, but this cannot happen by (a).
\end{proof}

\bl\label{ext4} Let $g$ be obtained from $f$ by a legal move of weight
zero. Then
$$\operatorname{Ext}^1(L(\lambda(g)),L(\lambda(f)))\neq 0.$$
\el
\begin{proof} Let the legal move have start $s_p$. By Lemma
  ~\ref{ext3} we have
$$\operatorname{Ext}^1_{{\sfq}^{(p-1)}}(L_{{\sfq}^{(p-1)}}(\lambda(g)),L_{{\sfq}^{(p-1)}}(\lambda(f)))\neq 0.$$
We will prove
$$\operatorname{Ext}^1_{{\sfq}^{(i)}}(L_{{\sfq}^{(i)}}(\lambda(g)),L_{{\sfq}^{(i)}}(\lambda(f)))\neq 0$$
for all $i\leq p-1$ by reverse induction in $i$. So we assume that the statement is
true for $i$ and prove it for $i-1$.
Consider a non-split exact sequence
$$0\to L_{{\sfq}^{(i)}}(\lambda(g))\to V \to L_{{\sfq}^{(i)}}(\lambda(f)) \to 0.$$
Apply $\Gamma^{(i)}$ to the sequence to obtain
$$0\to \Gamma^{(i)}(L_{{\sfq}^{(i)}}(\lambda(g)))\stackrel{\phi}{\to} \Gamma^{(i)}(V) \to \Gamma^{(i)}(L_{{\sfq}^{(i)}}(\lambda(f))) \to
0.$$
This sequence is not exact, but $\Gamma^{(i)}$ is exact on the
right, so we
have an isomorphism
$$\Gamma^{(i)}(V)/\operatorname{Im}\phi\simeq \Gamma^{(i)}(L_{{\sfq}^{(i)}}(\lambda(f))).$$
Now let $M$ and $N$ be the proper maximal submodules in  $\Gamma^{(i)}(L_{{\sfq}^{(i)}}(\lambda(g)))$
and $\Gamma^{(i)}(L_{{\sfq}^{(i)}}(\lambda(f)))$ respectively, and let
$X=\Gamma^{(i)}(V)/\phi(M)$. We have an exact sequence
$$0\to L_{{\sfq}^{(i-1)}}(\lambda(g)) \to X \stackrel{\pi}{\to} \Gamma^{(i)}(L_{{\sfq}^{(i)}}(\lambda(f)))\to 0.$$
From Theorem \ref{selecta2} we have that $$N= \bigoplus_{h\in \LM(f,i)}L_{{\sfq}^{(i-1)}}(\lambda(h)).$$
By Lemma \ref{extcomb} and Corollary \ref{corext1}
$$\operatorname{Ext}^1_{{\sfq}^{(i-1)}}(N,L_{{\sfq}^{(i-1)}}(\lambda(g)))=0.$$
Therefore
$$\pi^{-1}(N)=L_{{\sfq}^{(i-1)}}(\lambda(g))\oplus \bigoplus_{h\in \LM(f,i)}L_{{\sfq}^{(i-1)}}(\lambda(h)).$$
So $X/(\bigoplus_{h\in \LM(f,i)}L_{{\sfq}^{(i-1)}}(\lambda(h)))$ gives a
non-trivial extension between $L_{{\sfq}^{(i-1)}}(\lambda(g))$ and $L_{{\sfq}^{(i-1)}}(\lambda(f))$.
The case $i=0$ implies the statement.
\end{proof}

Corollary \ref{corext1}, Corollary \ref{corext2}, Lemma \ref{ext4} and
(\ref{duality}) imply Theorem B.\\
\br \label{rm69}{\rm  We thank Jon Brundan for pointing out that Theorem B can also be derived from
results in \cite{BS1}, \cite{BS2} and \cite{BS4}. The first two of these papers concern Khovanov's diagram algebra, and thus the remarks below rely
on Theorem 1.1 in \cite{BS4} which establishes
 an equivalence of categories between our category $\cF$ and the category of finite dimensional modules over a version of Khovanov's algebra studied in \cite{BS1}, \cite{BS2}.

It is shown in the proof of  \cite{BS2} Corollary 5.15, that the finite dimensional simple modules $\mathbb{L}(\gl)$  for Khovanov's algebra satisfy \be \label{Joneq}
\dim \Ext^1(\mathbb{L}(\gl),\mathbb{L}(\mu)) = p_{\gl,\mu}^{(1)}+p_{\mu,\gl}^{(1)}\ee
where  $p_{\gl,\mu}(q) = \sum_{i \geq 0}p_{\gl,\mu}^{(i)}q^i $ are the Kazhdan-Lusztig
polynomials considered there.

Hence assuming the equivalence of categories, Theorem B(a) follows from \cite{Br} Theorem 4.51.   Also from (\ref{Joneq}) and  the definition of a legal move of weight zero, Theorem B (b) follows from the assertion
\bl \label{jons} We have
$p_{\gl,\mu}^{(1)} \neq 0$ if and only if $\mu$ is obtained from $\gl$
by interchanging the labels at the ends of one of the caps in the cap
diagram of $\gl.$
\el
It is possible to prove  Lemma \ref{jons} using the recursive description
of $p_{\gl,\mu}^{(1)}$ from  \cite{LS} Lemma 6.6. The proof of this result is
reproduced in \cite{BS2}, Lemma 5.2.}\er

\end{document}